\documentclass[12pt]{article}%
\usepackage{amsfonts}
\usepackage{sw20bams}
\usepackage{amsmath}
\usepackage{amssymb}
\usepackage{graphicx}%
\providecommand{\U}[1]{\protect\rule{.1in}{.1in}}
\begin{document}

\title{Exercises in Iterational Asymptotics II}
\author{Steven Finch}
\date{March 6, 2025}
\maketitle

\begin{abstract}
The nonlinear recurrences we consider here include the functions $3x(1-x)$ and
$\cos(x)$, which possess attractive fixed points $2/3$ and $0.739...$
(Dottie's number). \ Detailed asymptotics for oscillatory convergence are
found, starting with a 1960 paper by Wolfgang Thron. \ Another function,
$x/\left(  1+x\ln(1+x)\right)  $, gives rise to a sequence with monotonic
convergence to $0$ but requires substantial work to calculate its associated
constant $C$.

\end{abstract}

\footnotetext{Copyright \copyright \ 2025 by Steven R. Finch. All rights
reserved.}

This paper is a continuation of \cite{F5-exc2}. \ As a preface, the quartic
recurrence%
\[%
\begin{array}
[c]{ccccc}%
x_{k}=x_{k-1}-a\,x_{k-1}^{3}+b\,x_{k-1}^{4}, &  & a>0, &  & b\neq0
\end{array}
\]
has asymptotic expansion%
\[
x_{k}\sim\frac{1}{\sqrt{2a}}\frac{1}{k^{1/2}}+\frac{b}{2a^{2}}\frac{1}%
{k}+\frac{-3a^{3}+2b^{2}}{8\sqrt{2}a^{7/2}}\frac{\ln(k)}{k^{3/2}}+\frac
{C}{k^{3/2}}%
\]
as $k\rightarrow\infty$, where $C$ is a constant that depends not only on $a$
\&\ $b$ but also on the initial value $x_{0}$. \ The first two coefficients
appeared in \cite{St-exc2, IS-exc2}, but also much earlier as a special case
of Theorem 5.1 in \cite{Thrn-exc2}. \ (Beware of a mistaken $k^{-3/2}$ second
order in \cite{IS-exc2}.) \ Proof of Thron's theorem \cite{Thrn-exc2} involves
what we call the brute-force matching-coefficient method \cite{Sc-exc2,
F1-exc2, F2-exc2, F3-exc2}. \ Using this on the quintic recurrence%
\[%
\begin{array}
[c]{ccccc}%
x_{k}=x_{k-1}-a\,x_{k-1}^{3}+b\,x_{k-1}^{4}+d\,x_{k-1}^{5}, &  & a>0, &  &
d\neq0
\end{array}
\]
we obtain a revised third-order term%
\[
\frac{-3a^{3}+2b^{2}+2a\,d}{8\sqrt{2}a^{7/2}}\,\frac{\ln(k)}{k^{3/2}}.
\]
Using this instead on the sextic recurrence%
\[%
\begin{array}
[c]{ccccc}%
x_{k}=x_{k-1}-a\,x_{k-1}^{3}+b\,x_{k-1}^{4}+d\,x_{k-1}^{5}+e\,x_{k-1}^{5}, &
& a>0, &  & e\neq0
\end{array}
\]
we obtain the same third-order term as before; more terms beyond $k^{-3/2}$
are possible:%
\[
\frac{-3a^{3}b+2b^{3}+2a\,b\,d}{8a^{5}}\,\frac{\ln(k)}{k^{2}}+\frac
{a^{3}b-3b^{3}+4\sqrt{2}a^{7/2}b\,c-3a\,b\,d-a^{2}e}{4a^{5}}\frac{1}{k^{2}}.
\]
We could go on, examining the septic analog and computing more terms. \ The
most interesting feature of these particular recurrences is the missing
$x_{k-1}^{2}$ term: the gap $2$ between first \&\ second exponents, followed
by uniform gaps $1$ thereafter, leads to analytical difficulties. \ The
Mavecha-Laohakosol algorithm \cite{dB-exc2, BR-exc2, ML-exc2}, employed
extensively in \cite{F4-exc2, F6-exc2}, does not apply here nor does it appear
easily generalizable. \ Hence the brute-force method is necessary in Sections
1 \&\ 2. \ 

\section{Quatri\`{e}me exercice}

\textbf{Describe in detail }the oscillatory convergence of%
\[%
\begin{array}
[c]{ccccc}%
x_{k}=3x_{k-1}(1-x_{k-1}) &  & \text{for }k\geq1\text{;} &  & x_{0}=\dfrac
{1}{2}%
\end{array}
\]
to its limiting value $2/3$. \ Determine numerically both%
\[
C_{\text{o}}=\lim_{k\rightarrow\infty}k^{3/2}\left[  \left(  x_{2k+1}-\frac
{2}{3}\right)  -\left(  \frac{1}{6}\frac{1}{k^{1/2}}-\frac{1}{24k}-\frac
{11}{192}\frac{\ln(k)}{k^{3/2}}\right)  \right]
\]
and%
\[
C_{\text{e}}=\lim_{k\rightarrow\infty}k^{3/2}\left[  \left(  \frac{2}%
{3}-x_{2k}\right)  -\left(  \frac{1}{6}\frac{1}{k^{1/2}}+\frac{1}{24k}%
-\frac{11}{192}\frac{\ln(k)}{k^{3/2}}\right)  \right]  .
\]
Using $C_{\text{o}}$ and $C_{\text{e}}$, find the asymptotic expansions of
$x_{2k+1}$ and of $x_{2k}$ to order $k^{-4}$.

Let $f(x)=3x(1-x)$. \ Note that%
\[
x_{0}=\frac{1}{2}<\frac{9}{16}=x_{2}<\ldots<\frac{2}{3}<\ldots<x_{3}%
=\dfrac{189}{256}<\dfrac{3}{4}=x_{1}%
\]
and%
\begin{align*}
f(f(x))  &  =9x(1-x)\left[  1-3x(1-x)\right] \\
&  =\left\{
\begin{array}
[c]{ccc}%
\dfrac{2}{3}+\left(  x-\dfrac{2}{3}\right)  -18\left(  x-\dfrac{2}{3}\right)
^{3}-27\left(  x-\dfrac{2}{3}\right)  ^{4} &  & \text{if }x>\dfrac{2}{3},\\
\dfrac{2}{3}-\left(  \dfrac{2}{3}-x\right)  +18\left(  \dfrac{2}{3}-x\right)
^{3}-27\left(  \dfrac{2}{3}-x\right)  ^{4} &  & \text{if }x<\dfrac{2}{3}.
\end{array}
\right.
\end{align*}
Setting $u_{k}=x_{2k+1}-\frac{2}{3}$ and $v_{k}=\frac{2}{3}-x_{2k}$, we obtain
recurrences
\[%
\begin{array}
[c]{ccc}%
u_{k}=u_{k-1}-18u_{k-1}^{3}-27u_{k-1}^{4}, &  & u_{0}=\dfrac{3}{4}-\dfrac
{2}{3}=\dfrac{1}{12},
\end{array}
\]%
\[
u_{k}\sim\frac{1}{6}\frac{1}{k^{1/2}}-\frac{1}{24k}-\frac{11}{192}\frac
{\ln(k)}{k^{3/2}}+\frac{C_{\text{o}}}{k^{3/2}}%
\]
and%
\[%
\begin{array}
[c]{ccc}%
v_{k}=v_{k-1}-18v_{k-1}^{3}+27v_{k-1}^{4}, &  & v_{0}=\dfrac{2}{3}-\dfrac
{1}{2}=\dfrac{1}{6},
\end{array}
\]%
\[
v_{k}\sim\frac{1}{6}\frac{1}{k^{1/2}}+\frac{1}{24k}-\frac{11}{192}\frac
{\ln(k)}{k^{3/2}}+\frac{C_{\text{e}}}{k^{3/2}}%
\]
with preliminary asymptotics from our preface.

As an aside, the short expansions for $u_{k}$ and $v_{k}$ are identical except
for a single sign ($-1/24$ for $u_{k}$ and $+1/24$ for $v_{k}$). \ We show
momentarily that longer expansions for $u_{k}$ and $v_{k}$ are similarly
identical -- corresponding coefficients are equal except possibly for sign --
we define the functions $x-18x^{3}-27x^{4}$ and $x-18x^{3}+27x^{4}$ to be
\textbf{kindred}. \ Many pairs of kindred functions were exhibited in
\cite{F4-exc2}, all resembling inverses, tailored vaguely. \ The inverse of
$x-18x^{3}-27x^{4}$, restricted to the interval $[0,1/12]$, is%
\[
\frac{1}{6}\left(  -1+\sqrt{3-2\sqrt{1-12x}}\right)  =x+18x^{3}+27x^{4}%
+972x^{5}+\cdots
\]
which, upon tailoring into an alternating series, gives $x-18x^{3}%
+27x^{4}-972x^{5}+\cdots$. \ It would seem that kindredness is more common
than once believed, and that a kindred triple (if not more)\ exists. \ End of aside.

To implement the brute-force matching-coefficient method for $u_{k}$, we
expand%
\[%
\begin{array}
[c]{ccccccccc}%
\dfrac{\ln(k+1)^{i}}{(k+1)^{2}}, &  & \dfrac{\ln(k+1)^{j}}{(k+1)^{5/2}}, &  &
\dfrac{\ln(k+1)^{\ell}}{(k+1)^{3}}, &  & \dfrac{\ln(k+1)^{m}}{(k+1)^{7/2}}, &
& \dfrac{\ln(k+1)^{n}}{(k+1)^{4}}%
\end{array}
\]
for $i=1,0$; $j=2,1,0$; $\ell=2,1,0$; $m=3,2,1,0$; $n=3,2,1,0$ and compare a
series for $u_{k+1}$ with a series for $u_{k}-18u_{k}^{3}-27u_{k}^{4}$. \ This
yields additional terms
\begin{align*}
&  \frac{11}{384}\frac{\ln(k)}{k^{2}}-\left(  \frac{5}{384}+\frac{C_{\text{o}%
}}{2}\right)  \frac{1}{k^{2}}+\frac{121}{4096}\frac{\ln(k)^{2}}{k^{5/2}%
}-\left(  \frac{121}{3072}+\frac{33C_{\text{o}}}{32}\right)  \frac{\ln
(k)}{k^{5/2}}\\
&  +\left(  \frac{77}{3072}+\frac{11C_{\text{o}}}{16}+9C_{\text{o}}%
^{2}\right)  \frac{1}{k^{5/2}}-\frac{121}{6144}\frac{\ln(k)^{2}}{k^{3}%
}+\left(  \frac{77}{2048}+\frac{11C_{\text{o}}}{16}\right)  \frac{\ln
(k)}{k^{3}}\\
&  -\left(  \frac{139}{6144}+\frac{21C_{\text{o}}}{32}+6C_{\text{o}}%
^{2}\right)  \frac{1}{k^{3}}-\frac{6655}{393216}\frac{\ln(k)^{3}}{k^{7/2}%
}+\left(  \frac{1331}{24576}+\frac{1815C_{\text{o}}}{2048}\right)  \frac
{\ln(k)^{2}}{k^{7/2}}\\
&  -\left(  \frac{2299}{32768}+\frac{121C_{\text{o}}}{64}+\frac{495C_{\text{o}%
}^{2}}{32}\right)  \frac{\ln(k)}{k^{7/2}}+\left(  \frac{2293}{73728}%
+\frac{627C_{\text{o}}}{512}+\frac{33C_{\text{o}}^{2}}{2}+90C_{\text{o}}%
^{3}\right)  \frac{1}{k^{7/2}}\\
&  +\frac{1331}{98304}\frac{\ln(k)^{3}}{k^{4}}-\left(  \frac{10285}%
{196608}+\frac{363C_{\text{o}}}{512}\right)  \frac{\ln(k)^{2}}{k^{4}}+\left(
\frac{297}{4096}+\frac{935C_{\text{o}}}{512}+\frac{99C_{\text{o}}^{2}}%
{8}\right)  \frac{\ln(k)}{k^{4}}\\
&  -\left(  \frac{9959}{294912}+\frac{81C_{\text{o}}}{64}+\frac{255C_{\text{o}%
}^{2}}{16}+72C_{\text{o}}^{3}\right)  \frac{1}{k^{4}}%
\end{align*}
in the expansion for $u_{k}$. \ 

Comparing likewise a series for $v_{k+1}$ with a series for $v_{k}-18v_{k}%
^{3}+27v_{k}^{4}$ yields additional terms%
\begin{align*}
&  -\frac{11}{384}\frac{\ln(k)}{k^{2}}+\left(  \frac{5}{384}+\frac
{C_{\text{e}}}{2}\right)  \frac{1}{k^{2}}+\frac{121}{4096}\frac{\ln(k)^{2}%
}{k^{5/2}}-\left(  \frac{121}{3072}+\frac{33C_{\text{e}}}{32}\right)
\frac{\ln(k)}{k^{5/2}}\\
&  +\left(  \frac{77}{3072}+\frac{11C_{\text{e}}}{16}+9C_{\text{e}}%
^{2}\right)  \frac{1}{k^{5/2}}+\frac{121}{6144}\frac{\ln(k)^{2}}{k^{3}%
}-\left(  \frac{77}{2048}+\frac{11C_{\text{e}}}{16}\right)  \frac{\ln
(k)}{k^{3}}\\
&  +\left(  \frac{139}{6144}+\frac{21C_{\text{e}}}{32}+6C_{\text{e}}%
^{2}\right)  \frac{1}{k^{3}}-\frac{6655}{393216}\frac{\ln(k)^{3}}{k^{7/2}%
}+\left(  \frac{1331}{24576}+\frac{1815C_{\text{e}}}{2048}\right)  \frac
{\ln(k)^{2}}{k^{7/2}}\\
&  -\left(  \frac{2299}{32768}+\frac{121C_{\text{e}}}{64}+\frac{495C_{\text{e}%
}^{2}}{32}\right)  \frac{\ln(k)}{k^{7/2}}+\left(  \frac{2293}{73728}%
+\frac{627C_{\text{e}}}{512}+\frac{33C_{\text{e}}^{2}}{2}+90C_{\text{e}}%
^{3}\right)  \frac{1}{k^{7/2}}\\
&  -\frac{1331}{98304}\frac{\ln(k)^{3}}{k^{4}}+\left(  \frac{10285}%
{196608}+\frac{363C_{\text{e}}}{512}\right)  \frac{\ln(k)^{2}}{k^{4}}-\left(
\frac{297}{4096}+\frac{935C_{\text{e}}}{512}+\frac{99C_{\text{e}}^{2}}%
{8}\right)  \frac{\ln(k)}{k^{4}}\\
&  +\left(  \frac{9959}{294912}+\frac{81C_{\text{e}}}{64}+\frac{255C_{\text{e}%
}^{2}}{16}+72C_{\text{e}}^{3}\right)  \frac{1}{k^{4}}%
\end{align*}
in the expansion for $v_{k}$. \ 

Our simple procedure for estimating the constant $C_{\text{o}}$ involves
computing $u_{K}$ exactly via recursion, for some suitably large index $K$.
\ We then set the value $u_{K}$ equal to our series and numerically solve for
$C_{\text{o}}$. \ When employing terms up to order $k^{-4}$, an index
$\approx10^{10}$ may be required for $25$ digits of precision in the
$C_{\text{o}}$ estimate:%
\[
C_{\text{o}}=-0.1805303007686495535981970....
\]
We likewise find $C_{\text{e}}$ from computing $v_{K}$:%
\[
C_{\text{e}}=-0.1388636341019828869315303....
\]
It is natural to speculate about the algebraic independence of these constants.

\section{Cinqui\`{e}me exercice}

\textbf{Describe in detail }the monotonic convergence of%
\[%
\begin{array}
[c]{ccccc}%
y_{k}=\dfrac{y_{k-1}}{1+y_{k-1}\ln(1+y_{k-1})} &  & \text{for }k\geq1\text{;}
&  & y_{0}=1
\end{array}
\]
to $0$. \ Determine numerically%
\[
C=\lim_{k\rightarrow\infty}k^{3/2}\left[  y_{k}-\left(  \frac{1}{\sqrt{2}%
}\frac{1}{k^{1/2}}+\frac{1}{4k}-\frac{7}{48\sqrt{2}}\frac{\ln(k)}{k^{3/2}%
}\right)  \right]
\]
and find the asymptotic expansion of $y_{k}$ to order $k^{-4}$.

A\ little background is helpful. \ Corollary 4 of \cite{P1-exc2} was devoted
to $x+\ln\left(  \alpha+\frac{\beta}{x}\right)  $ where $\alpha>1$ and
$\beta>0$. \ The case $\alpha=\beta=1$ is on the boundary of allowable values
and Popa's wide-ranging expansion for $x_{k}$ does not apply to $f(x)=x+\ln
\left(  1+\frac{1}{x}\right)  $. \ We focus on $y_{k}=1/x_{k}$, which
satisfies $y_{k}=g(y_{k-1})$ where%
\[
g(y)=\dfrac{1}{f\left(  1/y\right)  }=\dfrac{y}{1+y\ln(1+y)}=y-y^{3}+\dfrac
{1}{2}y^{4}+\dfrac{2}{3}y^{5}-\dfrac{3}{4}y^{6}-\dfrac{17}{60}y^{7}+\cdots.
\]
As before, the $y^{2}$ term is missing and hence Mavecha-Laohakosol is
inapplicable. \ In our prior work with brute force \cite{F1-exc2, F2-exc2,
F3-exc2}, we never once attempted analysis on a transcendental function.
\ Thus we study algebraic fits to $g(y)$ of polynomial degrees $4$, $5$, $6$,
$7$ and assess accuracy for various $y_{k}$-expansion lengths.

\subsection{Quartique}

As in Section 1, an estimate for the constant $C$ involves calculating $y_{K}$
exactly via recursion (logarithmic) for large $K$. \ The remaining steps are
to choose a Taylor approximation for $g(y)$ and then to select a cutoff for
the corresponding asymptotic series. \ Using the quartic $y-y^{3}+\tfrac{1}%
{2}y^{4}$, it is tempting to set the value $y_{K}$ equal to our series%
\[
\frac{1}{\sqrt{2}}\frac{1}{k^{1/2}}+\frac{1}{4k}-\frac{5^{\ast}}{16\sqrt{2}%
}\frac{\ln(k)}{k^{3/2}}+\frac{C}{k^{3/2}}%
\]
from the preface ($a=1$, $b=1/2$). \ This is unwise, however, as $C$ fails to
converge, seemingly increasing without bound. \ The coefficient $5/(16\sqrt
{2})$ is starred because it is only transient, i.e., based on $d=0$.

\subsection{Quintique}

Using the quintic $y-y^{3}+\tfrac{1}{2}y^{4}+\frac{2}{3}y^{5}$, we set the
value $y_{K}$ equal to our series%
\[
\frac{1}{\sqrt{2}}\frac{1}{k^{1/2}}+\frac{1}{4k}-\frac{7}{48\sqrt{2}}\frac
{\ln(k)}{k^{3/2}}+\frac{C}{k^{3/2}}%
\]
from the preface ($a=1$, $b=1/2$, $d=2/3$). \ This gives $C=-0.3318...$. \ If
we include the additional terms%
\[
-\frac{7}{96}\frac{\ln(k)}{k^{2}}+\left(  -\frac{7^{\ast}}{32}+\frac{C}%
{\sqrt{2}}\right)  \frac{1}{k^{2}}%
\]
then $C=-0.33181...$ emerges. \ The coefficient $7/32$ is only transient,
i.e., based on $e=0$.

\subsection{Sextique}

Using the sextic $y-y^{3}+\tfrac{1}{2}y^{4}+\frac{2}{3}y^{5}-\frac{3}{4}y^{6}%
$, we set the value $y_{K}$ equal to our series%
\[
\frac{1}{\sqrt{2}}\frac{1}{k^{1/2}}+\frac{1}{4k}-\frac{7}{48\sqrt{2}}\frac
{\ln(k)}{k^{3/2}}+\frac{C}{k^{3/2}}-\frac{7}{96}\frac{\ln(k)}{k^{2}}+\left(
-\frac{1}{32}+\frac{C}{\sqrt{2}}\right)  \frac{1}{k^{2}}%
\]
from the preface ($a=1$, $b=1/2$, $d=2/3$, $e=-3/4$). \ This gives
$C=-0.33181542...$. \ If we include the additional terms%
\begin{align*}
&  \frac{49}{1536\sqrt{2}}\frac{\ln(k)^{2}}{k^{5/2}}-\left(  \frac
{49}{1152\sqrt{2}}+\frac{7C}{16}\right)  \frac{\ln(k)}{k^{5/2}}+\left(
-\frac{43^{\ast}}{1152\sqrt{2}}+\frac{7C}{24}+\frac{3C^{2}}{\sqrt{2}}\right)
\frac{1}{k^{5/2}}\\
&  +\frac{49}{2304}\frac{\ln(k)^{2}}{k^{3}}-\left(  \frac{7}{2304}+\frac
{7C}{12\sqrt{2}}\right)  \frac{\ln(k)}{k^{3}}+\left(  -\frac{13^{\ast}}%
{2304}+\frac{C}{24\sqrt{2}}+2C^{2}\right)  \frac{1}{k^{3}}%
\end{align*}
then $C=-0.3318154296...$ emerges. \ The two starred coefficients are only transient.

\subsection{Septique}

Using the sextic $y-y^{3}+\tfrac{1}{2}y^{4}+\frac{2}{3}y^{5}-\frac{3}{4}%
y^{6}-\frac{17}{60}y^{7}$, we set the value $y_{K}$ equal to our series%
\begin{align*}
&  \frac{1}{\sqrt{2}}\frac{1}{k^{1/2}}+\frac{1}{4k}-\frac{7}{48\sqrt{2}}%
\frac{\ln(k)}{k^{3/2}}+\frac{C}{k^{3/2}}-\frac{7}{96}\frac{\ln(k)}{k^{2}%
}+\left(  -\frac{1}{32}+\frac{C}{\sqrt{2}}\right)  \frac{1}{k^{2}}\\
&  +\frac{49}{1536\sqrt{2}}\frac{\ln(k)^{2}}{k^{5/2}}-\left(  \frac
{49}{1152\sqrt{2}}+\frac{7C}{16}\right)  \frac{\ln(k)}{k^{5/2}}+\left(
-\frac{11}{5760\sqrt{2}}+\frac{7C}{24}+\frac{3C^{2}}{\sqrt{2}}\right)
\frac{1}{k^{5/2}}\\
&  +\frac{49}{2304}\frac{\ln(k)^{2}}{k^{3}}-\left(  \frac{7}{2304}+\frac
{7C}{12\sqrt{2}}\right)  \frac{\ln(k)}{k^{3}}+\left(  -\frac{3013}%
{11520}+\frac{C}{24\sqrt{2}}+2C^{2}\right)  \frac{1}{k^{3}}\\
&  -\frac{1715}{221184\sqrt{2}}\frac{\ln(k)^{3}}{k^{7/2}}+\left(  \frac
{343}{13824\sqrt{2}}+\frac{245C}{1536}\right)  \frac{\ln(k)^{2}}{k^{7/2}}\\
&  -\left(  \frac{203}{18432\sqrt{2}}+\frac{49C}{144}+\frac{35C^{2}}%
{16\sqrt{2}}\right)  \frac{\ln(k)}{k^{7/2}}+\left(  \frac{143}{3456\sqrt{2}%
}+\frac{29C}{384}+\frac{7C^{2}}{3\sqrt{2}}+5C^{3}\right)  \frac{1}{k^{7/2}}\\
&  -\frac{343}{55296}\frac{\ln(k)^{3}}{k^{4}}+\left(  \frac{833}{110592}%
+\frac{49C}{192\sqrt{2}}\right)  \frac{\ln(k)^{2}}{k^{4}}+\left(  \frac
{5453}{345600}-\frac{119C}{576\sqrt{2}}-\frac{7C^{2}}{4}\right)  \frac{\ln
(k)}{k^{4}}\\
&  -\left(  \frac{975007}{2304000}+\frac{779C}{3600\sqrt{2}}-\frac{17C^{2}%
}{24}-4\sqrt{2}C^{3}\right)  \frac{1}{k^{4}}%
\end{align*}
and $C=-0.331815429620156...$ emerges. \ This constant is unexpectedly
difficult to calculate:\ despite possessing the series to order $k^{-4}$, only
$15$ digits of $C$ are known.

We conclude that $C$ plays a role in the asymptotics of $x_{k}=1/y_{k}$ as
well:%
\[
x_{k}\sim\sqrt{2}k^{1/2}-\frac{1}{2}+\frac{7}{24\sqrt{2}}\frac{\ln(k)}%
{k^{1/2}}+\left(  \frac{1}{4\sqrt{2}}-2C\right)  \frac{1}{k^{1/2}}%
\]
but a general reciprocity formula (as in \cite{P1-exc2, P2-exc2} for a
specific scenario) seems out of reach. \ Also, for any integer $\ell\geq2$, a
gap $\ell+1$ between first \&\ second exponents in%
\[
\dfrac{y}{1+y^{\ell}\ln(1+y)}=y-y^{\ell+2}+\dfrac{1}{2}y^{\ell+3}-\dfrac{1}%
{3}y^{\ell+4}+\cdots
\]
opens the door to more related exploration.

\section{Sixi\`{e}me exercice}

\textbf{Consider} the famous recurrence%
\[%
\begin{array}
[c]{ccccc}%
x_{k}=\cos(x_{k-1}) &  & \text{for }k\geq1\text{;} &  & x_{0}=0.
\end{array}
\]
Quantify the convergence rate of $x_{k}$ as $k\rightarrow\infty$.

It is well known that%
\[
x_{0}=0<0.54\approx\cos(1)=x_{2}<\ldots<\theta<\ldots<x_{3}=\cos
(\cos(1))\approx0.85<1=x_{1}%
\]
where the limiting value%
\[
\theta=0.7390851332151606416553120...
\]
is Dottie's number \cite{Kapl-exc2, Pain-exc2}. \ Letting
\[%
\begin{array}
[c]{ccc}%
f(x)=\cos(\cos(\theta+x))-\theta, &  & g(x)=\theta-\cos(\cos(\theta-x))
\end{array}
\]
we have $x_{3}-\theta=f(x_{1}-\theta)$ and $\theta-x_{2}=g(\theta-x_{0})$.
\ The pattern is clear. \ Define%
\[%
\begin{array}
[c]{ccc}%
u_{k}=x_{2k+1}-\theta, &  & v_{k}=\theta-x_{2k}%
\end{array}
\]
and thus%
\[%
\begin{array}
[c]{ccc}%
u_{k+1}=f(u_{k}), &  & v_{k+1}=g(v_{k})
\end{array}
\]
for all $k$. \ Both $u_{k}$ and $v_{k}$ approach $0$; we determine the
respective speeds at which they do so, following Theorem 2.1 in
\cite{Thrn-exc2}. \ Note that $f(0)=g(0)=0$, $0<\max\{f(x),g(x)\}<x$ for all
$x>0$, and%
\[
f^{\prime}(0)=g^{\prime}(0)=1-\theta^{2}=0.4537531658603282480453425...<1.
\]
We now treat $f(x)$ and $g(x)$ separately.

The function%

\[
F(x)=\left\{
\begin{array}
[c]{ccc}%
\dfrac{f(x)-\left(  1-\theta^{2}\right)  x}{x^{2}} &  & \text{if }%
x>0,\medskip\\
\dfrac{\theta\sqrt{1-\theta^{2}}\left(  1-\sqrt{1-\theta^{2}}\right)  }{2} &
& \text{if }x=0
\end{array}
\right.
\]
is continuous and bounded on $[0,\infty)$; in fact, $\left\vert
F(x)\right\vert <M=0.27279$ by calculus. \ Observe that $\theta^{2}%
/(2M)\approx1.0012$ and hence $u_{k}<\theta^{2}/(2M)$ always. \ Because%
\[
\frac{u_{k+1}}{u_{k}}=\frac{f(u_{k})}{u_{k}}=\left(  1-\theta^{2}\right)
+F(u_{k})\,u_{k}<\left(  1-\theta^{2}\right)  +M\,\frac{\theta^{2}}%
{2M}=1-\frac{\theta^{2}}{2}%
\]
we have%
\[
u_{k+1}<\left(  1-\frac{\theta^{2}}{2}\right)  u_{k}<\left(  1-\frac
{\theta^{2}}{2}\right)  ^{2}u_{k-1}<\ldots<\left(  1-\frac{\theta^{2}}%
{2}\right)  ^{k+1}u_{0}.
\]
It follows that the series%
\[
\frac{1}{1-\theta^{2}}\,%
{\displaystyle\sum\limits_{k=0}^{\infty}}
\,u_{k}\left\vert F(u_{k})\right\vert <\frac{M}{1-\theta^{2}}\,%
{\displaystyle\sum\limits_{k=0}^{\infty}}
\,u_{k}<\frac{M\,u_{0}}{1-\theta^{2}}\,%
{\displaystyle\sum\limits_{k=0}^{\infty}}
\left(  1-\frac{\theta^{2}}{2}\right)  ^{k}%
\]
converges, which in turn implies that the product%
\[%
{\displaystyle\prod\limits_{k=0}^{\infty}}
\left(  1+\frac{1}{1-\theta^{2}}\,u_{k}\,F(u_{k})\right)
\]
also converges. \ Finally, multiplying both sides of
\[
\frac{1}{1-\theta^{2}}\frac{u_{j+1}}{u_{j}}=1+\frac{1}{1-\theta^{2}}%
\,u_{j}\,F(u_{j})
\]
from $j=0$ to $k-1$ gives%
\[
\frac{1}{(1-\theta^{2})^{k}}\,\frac{u_{k}}{u_{0}}=%
{\displaystyle\prod\limits_{j=0}^{k-1}}
\left(  1+\frac{1}{1-\theta^{2}}\,u_{j}\,F(u_{j})\right)
\]
and therefore%
\[
\lim_{k\rightarrow\infty}\frac{u_{k}}{(1-\theta^{2})^{k}}=(1-\theta)%
{\displaystyle\prod\limits_{j=0}^{\infty}}
\left(  1+\frac{1}{1-\theta^{2}}\,u_{j}\,F(u_{j})\right)
=0.2682998330950090571338993....
\]
Having finished with $f(x)$, we now investigate $g(x)$.

The function%

\[
G(x)=\left\{
\begin{array}
[c]{ccc}%
\dfrac{g(x)-\left(  1-\theta^{2}\right)  x}{x^{2}} &  & \text{if }%
x>0,\medskip\\
-\dfrac{\theta\sqrt{1-\theta^{2}}\left(  1-\sqrt{1-\theta^{2}}\right)  }{2} &
& \text{if }x=0
\end{array}
\right.
\]
is continuous and bounded on $[0,\infty)$; in fact, $\left\vert
G(x)\right\vert <M=0.30697$ by calculus. \ Observe that $\theta^{2}%
/(2M)\approx0.8897$ and hence $v_{k}<\theta^{2}/(2M)$ always. \ A similar line
of reasoning gives%
\[
\lim_{k\rightarrow\infty}\frac{v_{k}}{(1-\theta^{2})^{k}}=\theta\,%
{\displaystyle\prod\limits_{j=0}^{\infty}}
\left(  1+\frac{1}{1-\theta^{2}}\,v_{j}\,G(v_{j})\right)
=0.3983002403035094139563243....
\]
The two constants here differ by a factor of $\sqrt{1-\theta^{2}}$. \ 

\section{Septi\`{e}me exercice}

\textbf{Return to} the logistic map%
\[%
\begin{array}
[c]{ccccc}%
x_{k}=\lambda x_{k-1}(1-x_{k-1}) &  & \text{for }k\geq1\text{;} &  & 0<x_{0}<1
\end{array}
\]
where $1<\lambda<3$. \ Quantify the convergence rate of $x_{k}$ as
$k\rightarrow\infty$.

The limiting value $\mu=(\lambda-1)/\lambda$ satisfies $0<\mu<2/3$. \ We
initially examine $1<\lambda<2$. \ If $\ell(x)=\lambda\left(  x-x^{2}\right)
$, then $\ell(\mu)=\mu$ (being a fixed point),%
\[%
\begin{array}
[c]{ccccc}%
\ell^{\prime}(\mu)=\lambda(1-2\mu)=\lambda-2(\lambda-1)=2-\lambda, &  &
\ell^{\prime\prime}(\mu)/2=-\lambda, &  & \ell^{\prime\prime\prime}(\mu)/6=0
\end{array}
\]
and so%
\[
\ell(x)=\mu+(2-\lambda)(x-\mu)-\lambda(x-\mu)^{2}.
\]
Assume WLOG that $x_{0}>\mu$. \ The sequence $\{x_{k}\}$ is monotone
decreasing. \ Letting%
\[
f(x)=(2-\lambda)x-\lambda\,x^{2}%
\]
we have $x_{1}-\mu=f(x_{0}-\mu)$. \ Define $w_{0}=x_{0}-\mu$ and
$w_{k+1}=f(w_{k})$ for all $k$. \ The conditions for Theorem 2.1 in
\cite{Thrn-exc2} are met; in particular, $f^{\prime}(0)=2-\lambda<1$ and%
\[
F(x)=\dfrac{f(x)-\left(  2-\lambda\right)  x}{x^{2}}=-\lambda
\]
for all $x$. \ Convergence of the associated product follows as before. \ For
example, if $\lambda=3/2$ and $x_{0}=1/2$, then%
\[
\lim_{k\rightarrow\infty}\frac{w_{k}}{(2-\lambda)^{k}}=w_{0}\,%
{\displaystyle\prod\limits_{j=0}^{\infty}}
\left(  1-\frac{\lambda\,w_{j}}{2-\lambda}\right)
=0.0654844754592965980119173....
\]

The recurrence is trivial if $\lambda=2$:%
\[
x_{k}=\frac{1-\left(  1-2x_{0}\right)  ^{2^{k}}}{2}%
\]
as can be readily verified. \ Note the special cases $x_{0}=1/2$ and
$x_{0}=(1-e^{-1})/2$, for which%
\[%
\begin{array}
[c]{ccccc}%
x_{k}=\dfrac{1}{2}\text{ \ (identically)} &  & \text{and} &  & x_{k}=\dfrac
{1}{2}\left(  1-e^{-2^{k}}\right)
\end{array}
\]
respectively.

We finally examine $2<\lambda<3$. \ By the Chain Rule \cite{Jns-exc2,
HNY-exc2},%
\[
(\ell\circ\ell)^{\prime}(\mu)=\ell^{\prime}(\ell(\mu))\ell^{\prime}(\mu
)=\ell^{\prime}(\mu)^{2}=(2-\lambda)^{2},
\]%
\begin{align*}
(\ell\circ\ell)^{\prime\prime}(\mu)/2  &  =\left\{  \ell^{\prime\prime}%
(\ell(\mu))\ell^{\prime}(\mu)^{2}+\ell^{\prime}(\ell(\mu))\ell^{\prime\prime
}(\mu)\right\}  /2\\
&  =\ell^{\prime\prime}(\mu)\ell^{\prime}(\mu)\left[  \ell^{\prime}%
(\mu)+1\right]  /2\\
&  =(-2\lambda)(2-\lambda)\left(  3-\lambda\right)  /2,
\end{align*}%
\begin{align*}
(\ell\circ\ell)^{\prime\prime\prime}(\mu)/6  &  =\left\{  \ell^{\prime
\prime\prime}(\ell(\mu))\ell^{\prime}(\mu)^{3}+3\ell^{\prime\prime}(\ell
(\mu))\ell^{\prime}(\mu)\ell^{\prime\prime}(\mu)+\ell^{\prime}(\ell(\mu
))\ell^{\prime\prime\prime}(\mu)\right\}  /6\\
&  =\left\{  0+3\ell^{\prime\prime}(\mu)^{2}\ell^{\prime}(\mu)+0\right\}  /6\\
&  =3(-2\lambda)^{2}(2-\lambda)/6
\end{align*}
and so%
\[
\ell(\ell(x))=\mu+(\lambda-2)^{2}(x-\mu)-(\lambda-3)(\lambda-2)\lambda
(x-\mu)^{2}-2(\lambda-2)\lambda^{2}(x-\mu)^{3}-\lambda^{3}(x-\mu)^{4}.
\]
Assume WLOG that $x_{0}<\mu$. \ The sequence $\{x_{k}\}$ is oscillatory.
\ Letting%
\[
f(x)=(\lambda-2)^{2}x-(\lambda-3)(\lambda-2)\lambda\,x^{2}-2(\lambda
-2)\lambda^{2}x^{3}-\lambda^{3}x^{4},
\]%
\[
g(x)=(\lambda-2)^{2}x+(\lambda-3)(\lambda-2)\lambda\,x^{2}-2(\lambda
-2)\lambda^{2}x^{3}+\lambda^{3}x^{4}%
\]
we have $x_{3}-\mu=f(x_{1}-\mu)$ and $\mu-x_{2}=g(\mu-x_{0})$ \ Define
$u_{0}=x_{1}-\mu$,\ $v_{0}=\mu-x_{0}$ and $u_{k+1}=f(u_{k})$, $v_{k+1}%
=g(v_{k})$\ for all $k$. \ The conditions for Theorem 2.1 in \cite{Thrn-exc2}
are met; in particular, $f^{\prime}(0)=g^{\prime}(0)=(\lambda-2)^{2}<1$.
\ With%
\[%
\begin{array}
[c]{ccc}%
F(x)=\dfrac{f(x)-(\lambda-2)^{2}x}{x^{2}}, &  & G(x)=\dfrac{g(x)-(\lambda
-2)^{2}x}{x^{2}}%
\end{array}
\]
then taking $\lambda=5/2$ and $x_{0}=1/2$, we obtain convergent products
\[
\lim_{k\rightarrow\infty}\frac{u_{k}}{(2-\lambda)^{2k}}=u_{0}\,%
{\displaystyle\prod\limits_{j=0}^{\infty}}
\left(  1-\frac{1}{(2-\lambda)^{2}}u_{j}\,F(u_{j})\right)
=0.0266915553170954912963034...,
\]%
\[
\lim_{k\rightarrow\infty}\frac{v_{k}}{(2-\lambda)^{2k}}=v_{0}\,%
{\displaystyle\prod\limits_{j=0}^{\infty}}
\left(  1-\frac{1}{(2-\lambda)^{2}}v_{j}\,G(v_{j})\right)
=0.0533831106341909825926069....
\]
The constants here differ by a mere factor of $1/2$. \ This outcome is
completely unlike the mystery [surrounding iterates of $3x(1-x)$] that closes
Section 1.

\section{Acknowledgements}

The creators of Mathematica earn my gratitude every day:\ this paper could not
have otherwise been written.

\end{document}